\documentclass[3p]{elsarticle}
\usepackage{amsmath}
\usepackage{amssymb}
\usepackage{xspace}

\newtheorem{theorem}{Theorem}
\newtheorem{lemma}{Lemma}
\newtheorem{proposition}{Proposition}
\newtheorem{corollary}{Corollary}
\newdefinition{definition}{Definition}
\newdefinition{remark}{Remark}
\newdefinition{Example}{Example}
\newproof{proof}{Proof}
\newenvironment{Proof}{\begin{proof}}{\qed\end{proof}}

\newcommand{\field}{\mathbb{F}}
\newcommand{\two}{GF(2)}

\newcommand{\sdif}{\mathop{\mathrm{\Delta}}}

\begin{document}

\begin{frontmatter}

\title{The Nullity Theorem for Principal Pivot Transform}
\author{Robert Brijder}
\ead{robert.brijder@uhasselt.be}
\address{Hasselt University and Transnational University of Limburg, Belgium}

\begin{abstract}
We generalize the nullity theorem of Gustafson [Linear Algebra Appl.\ (1984)] from matrix inversion to principal pivot transform. Several special cases of the obtained result are known in the literature, such as a result concerning local complementation on graphs. As an application, we show that a particular matrix polynomial, the so-called nullity polynomial, is invariant under principal pivot transform.
\end{abstract}

\begin{keyword}
nullity theorem
\sep principal pivot transform
\sep Schur complement
\sep local complementation
\sep nullity polynomial
\vspace{0.3cm}
\MSC 05C50 \sep 15A03 \sep 15A09
\end{keyword}

\end{frontmatter}

\section{Introduction}
Motivated by the well-known linear complementarity problem, Tucker \cite{tucker1960} defined a matrix operation to study combinatorial equivalence of matrices. A slight modification of this matrix operation (some signs are different in the definition) became known as \emph{principal pivot transform}. Principal pivot transform partially (component-wise) inverts a given matrix along a given set of indices, and it is applied in various settings such as mathematical programming and numerical analysis, see \cite{Tsatsomeros2000151} for an overview.

The nullity theorem \cite{Gustafson1984}, independently discovered in \cite[Theorem~2]{Fiedler1986}, establishes a one-to-one correspondence between the submatrices of a nonsingular matrix and the submatrices of its inverse such that the nullities of the submatrices are retained by the correspondence. The power of the nullity theorem is well illustrated in \cite{Strang2004}. The main result of this paper (Theorem \ref{thm:null_thm_ppt}) generalizes the nullity theorem to principal pivot transform (for which matrix inversion is a special case). We show that several other special cases of this main result are known in the literature, including a result on graphs, and we show that a particular matrix polynomial, the so-called nullity polynomial, is invariant under principal pivot transform (Corollary~\ref{cor:null_pol_invariant}).

\section{Notation and Terminology}
For finite sets $U$ and $V$, a $U \times V$-matrix $A$ (over
some field $\field$) is a matrix where the rows are indexed by
$U$ and the columns are indexed by $V$, i.e., $A$ is formally a function $U
\times V \rightarrow \field$. Hence, the order of the
rows/columns is not fixed (i.e., interchanging rows or columns
is mute). Note that, e.g., the rank $r(A)$, the nullity (i.e.,
dimension of the null space) $n(A)$, and the inverse $A^{-1}$
are well defined for $A$ (the latter, of course, only when $A$
is square and nonsingular). We denote for $i \in U$ and $j \in V$,
the value of the $(i,j)^{th}$ entry of $A$ by $A[i,j]$. For $X \subseteq U$ and $Y \subseteq V$, the
submatrix of $A$ induced by $(X,Y)$, denoted by $A[X,Y]$, is
the restriction of $A$ to $X \times Y$.

Similarly, a vector indexed by $V$ is formally a function $V
\rightarrow \field$, and we denote the element of $v$
corresponding to $i \in V$ by $v[i]$. As usual, the family of
vectors indexed by $V$ is denoted by $\field^V$. For $Y
\subseteq V$, we let $\iota_{Y,V}$ be the usual injection
$\field^Y \rightarrow \field^V$ by padding zeros. More
precisely, $(\iota_{Y,V}(w))[x] = w[x]$ if $x \in Y$ and
$(\iota_{Y,V}(w))[x] = 0$ if $x \in V \setminus Y$. Similarly,
we let $\pi_{Y,V}$ be the usual projection of $\field^V
\rightarrow \field^Y$ by disregarding the entries with indices
in $V \setminus Y$. If $V$ is clear from the context, we simply
write $\iota_{Y}$ and $\pi_{Y}$ for $\iota_{Y,V}$ and
$\pi_{Y,V}$, respectively.

\section{Principal Pivot Transform} \label{sec:def_pivots}
In this section we recall principal pivot transform, which is
an operation for square matrices, see \cite{Tsatsomeros2000151}
for an overview.

Let $A$ be a $V \times V$-matrix (over an arbitrary field
$\field$) and let $X \subseteq V$ be such that the
corresponding principal submatrix $A[X,X]$ is nonsingular. The
\emph{principal pivot transform} (\emph{PPT}
for short) of $A$ on $X$, denoted by $A*X$, is defined as
follows:
\begin{eqnarray}\label{pivot_def}
A*X = \bordermatrix{
& \scriptstyle X & \scriptstyle V\setminus X \cr
\scriptstyle X & A[X,X]^{-1} & -A[X,X]^{-1} A[X,V\setminus X] \cr
\scriptstyle V\setminus X & A[V\setminus X,X] A[X,X]^{-1} & A[V\setminus X,V\setminus X] - A[V\setminus X,X] A[X,X]^{-1} A[X,V\setminus X]
}.
\end{eqnarray}
Matrix $A*X[V\setminus X,V\setminus X]$ is called the \emph{Schur complement}
of $A[X,X]$ in $A$ \cite{SchurBook2005}.

Principal pivot transform can be considered a partial inverse, as $A$ and $A*X$
are related as follows, where the vectors $x_1$
and $y_1$ correspond to the elements of $X$:
\begin{eqnarray} \label{pivot_def_reverse} A \left(
\begin{array}{c} x_1 \\ x_2 \end{array} \right) = \left(\begin{array}{c} y_1 \\ y_2 \end{array} \right) \mbox{ if and only if } A*X \left(
\begin{array}{c} y_1 \\ x_2 \end{array} \right) = \left(\begin{array}{c} x_1 \\ y_2 \end{array}
\right). \end{eqnarray}
Equation~(\ref{pivot_def_reverse}) characterizes PPT, see \cite[Theorem~3.1]{Tsatsomeros2000151}. Note
that if $A$ is nonsingular, then $A * V = A^{-1}$. Also note that by
Equation~(\ref{pivot_def_reverse}) PPT is
an involution (operation of order $2$), and more generally, if
$(A*X)*Y$ is defined, then it is equal to $A*(X \sdif Y)$, where $\sdif$ denotes symmetric difference.

We denote by $A \sharp X$ the matrix obtained from $A$ by replacing
every row of $A$ with index $x \in V \setminus X$ by
$i_x^T$ where $i_x$ is the vector having value $1$ at index $x$
and $0$ elsewhere. Note that $A \sharp X\left(\begin{array}{c} x_1 \\ x_2 \end{array} \right) = \left(\begin{array}{c} y_1 \\ x_2 \end{array} \right)$
with $x_1$, $x_2$, and $y_1$ from (\ref{pivot_def_reverse}). From this it follows that if $A[X,X]$ is nonsingular, then $(A \sharp X)^{-1} = (A*X)\sharp X$.

\section{Nullity Theorem for Principal Pivot Transform} \label{sec:null_thm_ppt}

The following theorem is used in \cite{BH/PivotNullityInvar/09} to generalize the recursive relation for interlace polynomials from graphs \cite{Arratia2004199,Aigner200411} to arbitrary square matrices over arbitrary fields.
\begin{proposition}[\cite{BH/PivotNullityInvar/09}]\label{prop:null_thm_ppt_psub}
Let $A$ be a $V\times V$-matrix (over some field) and let $Z
\subseteq V$ be such that $A[Z,Z]$ is nonsingular. Then, for
all $X \subseteq V$, $n(A*Z[X,X]) = n(A[X \sdif Z,X \sdif
Z])$.
\end{proposition}

We now recall the nullity theorem.
\begin{proposition}[The nullity theorem, \cite{Gustafson1984}]\label{prop:null_thm_inv}
Let $A$ be a nonsingular $V\times V$-matrix (over some field). Then, for all $X,Y \subseteq V$,
$n(A^{-1}[X,Y]) = n(A[V \setminus Y,V \setminus X])$.
\end{proposition}

\begin{remark}
We remark that the nullity theorem is in the literature often
formulated in an unnecessary cumbersome way, since rows and columns
of matrices are conventionally ordered and moreover not
explicitly indexed. Instead, the concise description of the
nullity theorem in Proposition~\ref{prop:null_thm_inv} shows
its nature as an elementary result.
\end{remark}

Since $A * V = A^{-1}$, it is natural to ask whether
Propositions~\ref{prop:null_thm_ppt_psub} and
\ref{prop:null_thm_inv} may be put under a common umbrella. We
now show that this is indeed the case. Moreover, we relate the null spaces of the submatrices of $A$ and $A*Z$.
\begin{theorem}[Nullity theorem for PPT]\label{thm:null_thm_ppt}
Let $A$ be a $V\times V$-matrix (over some field) and let $Z
\subseteq V$ be such that $A[Z,Z]$ is nonsingular. Then, for
all $X,Y \subseteq V$,
\begin{eqnarray}
\ker(A*Z[X,Y]) &=& \ker(A[X \sdif R, Y \sdif R] \cdot (A\sharp Z)^{-1}[Y \sdif R, Y]) \mbox{ and} \label{eqn_ker_thm}\\
n(A*Z[X,Y]) &=& n(A[X \sdif R, Y \sdif
R]), \label{eqn_null_thm}
\end{eqnarray}
where $R = Z \setminus (X \sdif Y)$.
\end{theorem}
\begin{Proof}
Let $v \in \ker(A[X,Y])$. Then $A[X,Y]v = 0$, and so $\pi_X(A\iota_Y(v))=0$.

If $Ax = y$, then for all $i \in V$,
\begin{eqnarray} \label{eqn_sharp}
((A \sharp Z)x)[i] = \begin{cases} y[i] & \mbox{if } i \in Z,\\
x[i] & \mbox{otherwise},
\end{cases}
\end{eqnarray}
and, by Equation~(\ref{pivot_def_reverse}),
\begin{eqnarray} \label{eqn_pivot_sharp}
((A*Z)(A \sharp Z)x)[i] = \begin{cases} x[i] & \mbox{if } i \in Z,\\
y[i] & \mbox{otherwise}.
\end{cases}
\end{eqnarray}

By Equation~(\ref{eqn_sharp}), $((A \sharp Z)\iota_Y(v))[i] = 0$ if $i \in (X \cap Z) \cup ((V\setminus Y) \setminus Z) = V \setminus (Y \sdif R)$. Thus $(A*Z \cdot A \sharp Z)\iota_Y(v) = (A*Z[V,Y\sdif R] \cdot A \sharp Z [Y\sdif R,Y])v$.

By Equation~(\ref{eqn_pivot_sharp}), $((A*Z \cdot A \sharp Z)\iota_Y(v))[i] = 0$ if $i \in (X\setminus Z)\cup ((V\setminus Y)\cap Z) = X \sdif R$. Thus $(A*Z[X\sdif R,Y\sdif R] \cdot A \sharp Z [Y\sdif R,Y])v = \pi_{X\sdif R}(A*Z \cdot A \sharp Z\iota_Y(v)) = 0$. Consequently, $v \in \ker(A*Z[X \sdif R, Y \sdif R] \cdot A\sharp Z[Y \sdif R, Y])$. Hence
\begin{eqnarray} \label{eqn_first_incl}
\ker(A[X,Y]) \subseteq \ker(A*Z[X \sdif R, Y \sdif R] \cdot A\sharp Z[Y \sdif R, Y]).
\end{eqnarray}
We show that if $v_1, v_2 \in \ker(A[X,Y])$ are distinct, then $A\sharp Z[Y \sdif R, Y]v_1 \neq A\sharp Z[Y \sdif R, Y]v_2$. Assume to the contrary that $A\sharp Z[Y \sdif R, Y]v_1 = \allowbreak A\sharp Z \allowbreak[Y \sdif R, Y]v_2$. Recall that $((A \sharp Z)\iota_Y(v))[i] = 0$ if $i \in V \setminus (Y \sdif R)$. Hence, $A\sharp Z \iota_Y(v_1) = A\sharp Z \iota_Y(v_2)$, and thus $\iota_Y(v_1) = \iota_Y(v_2)$ as $A\sharp Z$ is nonsingular. Consequently, $v_1 = v_2$ --- a contradiction. We thus obtain, by (\ref{eqn_first_incl}), $n(A*Z[X \sdif R, Y \sdif R]) \geq \allowbreak  n(A[X,Y])$.

We apply now (\ref{eqn_first_incl}) to $A:=A*Z$, $X:=X \sdif R$, and $Y:=Y \sdif R$. We obtain $\ker(A*Z[X \sdif R, Y \sdif R]) \subseteq \ker((A*Z)*Z[(X \sdif R)\allowbreak \sdif R', (Y \sdif R) \sdif R'] \cdot A*Z\sharp Z[(Y \sdif R) \sdif R', Y \sdif R])$ where $R' = Z \setminus ((X\sdif R) \sdif \allowbreak(Y\sdif R)) = Z \setminus (X \sdif Y) = R$. Hence we have
\begin{eqnarray} \label{eqn_second_incl}
\ker(A*Z[X \sdif R, Y \sdif R]) \subseteq \ker(A[X,Y] \cdot (A\sharp Z)^{-1}[Y, Y  \sdif R])
\end{eqnarray}
as $(A\sharp Z)^{-1} = A*Z\sharp Z$. We again obtain that if $v_1, v_2 \in \ker(A*Z[X \sdif R, Y \sdif R])$ are distinct, then $(A\sharp Z)^{-1}[Y, Y  \sdif R]v_1 \neq (A\sharp Z)^{-1}[Y, Y  \sdif R]v_2$. Thus, by (\ref{eqn_second_incl}), $n(A[X,Y]) \geq n(A*Z[X \sdif R, Y \sdif R])$.

Consequently, $n(A[X,Y]) = n(A*Z[X \sdif R, Y \sdif
R])$ and the inclusions of (\ref{eqn_first_incl}) and (\ref{eqn_second_incl}) are equalities. By change of variables $A:=A*Z$, we have (\ref{eqn_ker_thm}) and (\ref{eqn_null_thm}).
\end{Proof}
Note that Propositions~\ref{prop:null_thm_ppt_psub} and
\ref{prop:null_thm_inv} are Equation~(\ref{eqn_null_thm}) of Theorem~\ref{thm:null_thm_ppt} for
the cases $X = Y$ and $Z = V$, respectively. The case $Y = V
\setminus X$ is also of particular interest, and so we
explicitly state it here.

\begin{corollary}\label{cor:null_thm_ppt}
Let $A$ be a $V\times V$-matrix (over some field) and let $Z
\subseteq V$ be such that $A[Z,Z]$ is nonsingular. Then, for
all $X \subseteq V$, $n(A*Z[X,V \setminus X]) = n(A[X,V
\setminus X])$.
\end{corollary}
Equivalently, we have $r(A*Z[X,V \setminus X]) = r(A[X,V
\setminus X])$.

Special cases of Corollary~\ref{cor:null_thm_ppt} have been
considered in the literature.
Oum~\cite[Corollary~4.14]{Sangil/LAA/2011/QuasiOrdering} shows,
through the use of Lagrangian chain-groups, that Corollary~\ref{cor:null_thm_ppt} holds
for the case where
$A$ is skew-symmetric or symmetric. In the next section we show
that a result on graphs of Bouchet
\cite{Bouchet/SIAM/1987/DigraphDecompEulerSys} can also be seen
as a special case of Corollary~\ref{cor:null_thm_ppt}.

The polynomial $q(A) = \sum_{X \subseteq V} y^{n(A[X,X])}$ is a straightforward generalization of the interlace polynomial \cite{Arratia2004199,Aigner200411} and is shown to be invariant under PPT \cite{BH/PivotNullityInvar/09}. Due to Theorem~\ref{thm:null_thm_ppt}, we may now define another polynomial that is invariant under PPT. Let us define the (extended)
\emph{nullity polynomial} for a $V\times V$-matrix $A$ by
$$
p(A) = \sum_{X,Y \subseteq V} y^{n(A[X,Y])}.
$$
\begin{lemma} \label{lem:one_to_one}
For each $Z \subseteq V$, the function $f_Z: 2^V \times 2^V \rightarrow 2^V \times 2^V$ defined by $f_Z(X,Y) = (X \sdif R_{X,Y,Z}, \allowbreak Y \sdif
R_{X,Y,Z})$ with $R_{X,Y,Z} = Z \setminus (X \sdif Y)$ is an one-to-one correspondence.
\end{lemma}
\begin{Proof}
Since the domain and codomain of $f_Z$ are finite and equal, it suffices to show that $f_Z$ is injective. Assume that $f_Z(X,Y) = f_Z(X',Y')$. Then $X \sdif R_{X,Y,Z} = X' \sdif R_{X',Y',Z}$ and $Y \sdif R_{X,Y,Z} = Y' \sdif R_{X',Y',Z}$. Thus $X \sdif Y = (X \sdif R_{X,Y,Z}) \sdif (Y \sdif R_{X,Y,Z}) = (X' \sdif R_{X',Y',Z}) \sdif (Y' \sdif R_{X',Y',Z}) = X' \sdif Y'$. Therefore $R_{X,Y,Z} = R_{X',Y',Z}$, and so $X = X'$ and $Y = Y'$.
\end{Proof}

By Lemma~\ref{lem:one_to_one} and Theorem~\ref{thm:null_thm_ppt}, for each $i \in
\{0,\ldots,|V|\}$, the number of submatrices of $A$ of nullity
$i$ is invariant under PPT. Hence we have the following.

\begin{corollary} \label{cor:null_pol_invariant}
Let $A$ be a $V\times V$-matrix, and $Z \subseteq V$ such that $A[Z,Z]$ is nonsingular. Then $p(A) = p(A*Z)$.
\end{corollary}

\section{Graphs}
Let $G = (V,E)$ be a simple graph, i.e., without loops or
parallel edges. We write $V(G) = V$ and $E(G) = E$. The
\emph{neighborhood} of $v \in V$ in $G$, denoted by $N_G(v)$, is
$\{w \in V \mid \{v,w\} \in E(G)\}$. The \emph{local
complement} of $G$ at $v \in V$, denoted by $G^v$, is obtained
from $G$ by replacing the subgraph of $G$ induced by $N_G(v)$
by its complementary subgraph. Hence, if $u,w \in N_G(v)$ are distinct, then
$\{u,w\} \in E(G)$ if and only if $\{u,w\} \not\in E(G^v)$. Graphs $G$ and
$G'$ are said to be \emph{locally equivalent} if there is a
(possibly empty) sequence of local complementations such that
$G'$ is obtained from $G$. Since local complementation is an
involution, local equivalence induces an equivalence relation.

The \emph{adjacency matrix} $A(G)$ of $G$ is the $V(G) \times
V(G)$-matrix over $\two$ where for all $u,v \in V$, $A(G)[u,v]
= 1$ if and only if $\{u,v\} \in E(G)$ (note that the diagonal entries are
$0$). The following result is from Bouchet
\cite{Bouchet/SIAM/1987/DigraphDecompEulerSys} (see also
\cite[Section~3]{Bouchet/DM/1993/RecognLocEquivGraphs}), and is rediscovered in \cite[Proposition~2.6]{DBLP:journals/jct/Oum05}.

\begin{proposition}[\cite{Bouchet/SIAM/1987/DigraphDecompEulerSys}]\label{prop:bouchet_lc}
Let $G$ and $G'$ be simple graphs that are locally equivalent.
Then, for all $X \subseteq V$, $r(A(G')[X,V(G) \setminus X]) =
r(A(G)[X,V(G) \setminus X])$.
\end{proposition}

We remark that the function which, for a simple graph $G$, assigns every $X \subseteq V(G)$ to the value $r(A(G)[X,\allowbreak V(G) \setminus X])$, is called the \emph{connectivity function} in \cite{Bouchet/SIAM/1987/DigraphDecompEulerSys} and the \emph{cut-rank} in \cite{DBLP:journals/jct/Oum05}.

Let, for $X \subseteq V$, $I_X$ be the $V(G) \times
V(G)$-matrix over $\two$ where for
all $u,v \in V(G)$, $I_X[u,v] = 1$ if and only if $u = v \in X$. By Equation~(\ref{pivot_def}) it is easy to see that $A(G^v) = ((A(G)+I_{\{v\}})*\{v\})+I_{N_G(v)\cup\{v\}}$ for all $v \in V(G)$ (see also, e.g., \cite{BH/PivotLoopCompl/09}).
Proposition~\ref{prop:bouchet_lc} follows now readily from
Corollary~\ref{cor:null_thm_ppt}.

\begin{Proof}
[of Proposition~\ref{prop:bouchet_lc}] It suffices to consider the case $G' = G^v$ with $v \in V(G)$ as the general case follows by iteration. We have $A(G^v)[X,V(G) \setminus X] =
(((A(G)+I_{\{v\}})*\{v\})+I_{N_G(v)\cup\{v\}})[X,V(G) \setminus X] =
((A(G)+I_{\{v\}})*\{v\})[X,V(G) \setminus X]$. By
Corollary~\ref{cor:null_thm_ppt}, $r(((A(G)+I_{\{v\}})*\{v\})[X,V(G)
\setminus X]) = r((A(G)+I_{\{v\}})[X,V(G) \setminus X]) =
r(A(G)[X,V(G) \setminus X])$.
\end{Proof}

\subsection*{Acknowledgements}
We thank Hendrik Jan Hoogeboom for careful reading and useful comments on the paper.

\bibliography{nullity_thm_ppt}

\begin{thebibliography}{10}

\bibitem{Aigner200411}
M.~Aigner and H.~van~der Holst.
\newblock Interlace polynomials.
\newblock {\em Linear Algebra and its Applications}, 377:11--30, 2004.

\bibitem{Arratia2004199}
R.~Arratia, B.~Bollob\'as, and G.B. Sorkin.
\newblock The interlace polynomial of a graph.
\newblock {\em Journal of Combinatorial Theory, Series B}, 92(2):199--233,
  2004.

\bibitem{Bouchet/SIAM/1987/DigraphDecompEulerSys}
A.~Bouchet.
\newblock Digraph decompositions and eulerian systems.
\newblock {\em SIAM Journal on Algebraic and Discrete Methods}, 8(3):323--337,
  1987.

\bibitem{Bouchet/DM/1993/RecognLocEquivGraphs}
A.~Bouchet.
\newblock Recognizing locally equivalent graphs.
\newblock {\em Discrete Mathematics}, 114:75--86, 1993.

\bibitem{BH/PivotLoopCompl/09}
R.~Brijder and H.J. Hoogeboom.
\newblock The group structure of pivot and loop complementation on graphs and
  set systems.
\newblock {\em European Journal of Combinatorics}, 32:1353--1367, 2011.

\bibitem{BH/PivotNullityInvar/09}
R.~Brijder and H.J. Hoogeboom.
\newblock Nullity invariance for pivot and the interlace polynomial.
\newblock {\em Linear Algebra and its Applications}, 435:277--288, 2011.

\bibitem{Fiedler1986}
M.~Fiedler and T.L. Markham.
\newblock Completing a matrix when certain entries of its inverse are
  specified.
\newblock {\em Linear Algebra and its Applications}, 74:225--237, 1986.

\bibitem{Gustafson1984}
W.H. Gustafson.
\newblock A note on matrix inversion.
\newblock {\em Linear Algebra and its Applications}, 57:71--73, 1984.

\bibitem{DBLP:journals/jct/Oum05}
S.~Oum.
\newblock Rank-width and vertex-minors.
\newblock {\em Journal of Combinatorial Theory, Series B}, 95(1):79--100, 2005.

\bibitem{Sangil/LAA/2011/QuasiOrdering}
S.~Oum.
\newblock Rank-width and well-quasi-ordering of skew-symmetric or symmetric
  matrices.
\newblock {\em Linear Algebra and its Applications}, 436(7):2008--2036, 2012.

\bibitem{Strang2004}
G.~Strang and T.~Nguyen.
\newblock The interplay of ranks of submatrices.
\newblock {\em SIAM Review}, 46(4):637--646, 2004.

\bibitem{Tsatsomeros2000151}
M.J. Tsatsomeros.
\newblock Principal pivot transforms: properties and applications.
\newblock {\em Linear Algebra and its Applications}, 307(1-3):151--165, 2000.

\bibitem{tucker1960}
A.W. Tucker.
\newblock A combinatorial equivalence of matrices.
\newblock In {\em Combinatorial Analysis, Proceedings of Symposia in Applied
  Mathematics}, volume~X, pages 129--140. American Mathematical Society, 1960.

\bibitem{SchurBook2005}
F.~Zhang.
\newblock {\em The Schur Complement and Its Applications}.
\newblock Springer, 2005.

\end{thebibliography}

\end{document}